\documentclass[journal]{IEEEtran}
\IEEEoverridecommandlockouts
\usepackage{amsmath,amssymb,latexsym,epsfig,psfrag,graphicx}
\newtheorem{prp}{Proposition}
\newtheorem{thm}[prp]{Theorem}

\newenvironment{pf}{\smallbreak\noindent{\it Proof. }}{\hfill$\Box$\smallbreak}
\newenvironment{pf*}[1]{\smallbreak\noindent{\it #1}}{\hfill$\Box$\smallbreak}
\newcounter{definition}

\newcounter{remark}
\newenvironment{rmk}{\addtocounter{remark}{1}\smallbreak\noindent
  {\em Remark \theremark.}}{\smallbreak}
\newcounter{example}
\newenvironment{ex}{\addtocounter{example}{1}\smallbreak\noindent
  {\bf Example \theexample.}}{\hfill$\Box$\smallbreak}
\newenvironment{ex*}[1]{\addtocounter{example}{1}\smallbreak\noindent
  {\bf Example \theexample. {\bf #1}}}{\hfill$\Box$\smallbreak}

\newcommand{\realR}{\mathbb{R}}
\newcommand{\complexC}{\mathbb{C}}

\newcommand{\RHinfty}{\mathbb{RH}_{\infty}}
\newcommand{\DHinfty}{\mathbb{DH}_{\infty}}

\DeclareMathOperator{\diag}{diag}
\DeclareMathOperator{\ind}{ind}

\DeclareMathOperator{\trace}{trace }
  \title{Distributed Control of Positive Systems}
  \author{Anders Rantzer
  \thanks{A. Rantzer is with Automatic Control LTH,
    Lund University, Box 118, SE-221 00 Lund, Sweden. 
    Email: {\footnotesize\tt rantzer@control.lth.se}.
}}
\begin{document}
\maketitle
\begin{abstract}
A system is called positive if the set of non-negative states is left invariant by the dynamics.  
Stability analysis and controller optimization are greatly
simplified for such systems. For example, linear Lyapunov functions
and storage functions can be used instead of quadratic ones. This
paper shows how such methods can be used for synthesis of
distributed controllers. It also shows that stability and
performance of such control systems can be verified with a complexity that scales linearly with the
number of interconnections. Several results regarding scalable
synthesis and verfication are derived, including a new
stronger version of the Kalman-Yakubovich-Popov lemma for positive systems.
Some main results are stated for frequency domain models using the notion of
positively dominated system. The analysis is illustrated with applications to transportation
networks, vehicle formations and power systems. 
\end{abstract}

\section{Introduction}
Classical methods for multi-variable control, such as LQG and
$H_{\infty}$-optimization, suffer from a lack of scalability that make them
hard to use for large-scale systems.
The difficulties are partly due to computational
complexity, partly absence of distributed structure in the resulting
controllers. The complexity growth can be traced back to the fact that
stability verification of a linear system  
with $n$ states generally requires a Lyapunov function 
involving $n^2$ quadratic terms, even if the system
matrices are
sparse. The situation improves drastically if we restrict
attention to closed loop dynamics described by system matrices with nonnegative 
off-diagonal entries. Then stability and performance can be
verified using a Lyapunov function with only $n$ linear terms.
Sparsity can be exploited in performance verification and even
synthesis of distributed controllers can be done with a
complexity that grows linearly with the number of nonzero entries in
the system matrices. These observations have far-reaching
implications for control engineering:
\begin{enumerate}
\item The conditions that enable scalable solutions hold naturally in many important
  application areas, such as stochastic systems, economics,
  transportation networks, chemical reactions, power systems and ecology.
\item The essential mathematical property can be extended to frequency
  domain models. A sufficient
  condition is that the transfer functions involved are ``positively
  dominated''.  
\item In control applications, the assumption of positive dominance
  need not hold for the open loop process. However, a large-scale
  control system can often be 
  structured into local control loops that give positive
  dominance, thus enabling scalable methods for optimization of
  the global performance.
\end{enumerate}
The study of matrices with nonnegative coefficients has a long history,
dating back to the Perron-Frobenius Theorem in 1912. A classic book on
the topic is \cite{Berman+94}. The theory is used in Leontief
economics \cite{Leontief86}, where the states denote nonnegative quantities of 
commodities. It appears in the study of Markov chains
\cite{Seneta81}, where the states denote nonnegative probabilities and
in compartment models \cite{Jacquez96}, where the states could denote
populations of species. A nice introduction to
the subject is given in \cite{Luenberger79}. 
characterized by the property that a partial ordering of initial
states is preserved by the dynamics. Such dynamical systems were studied in a series
of papers by Hirsch, for example showing that
monotonicity generally implies convergence
almost everywhere \cite{Hirsch83,Hirsch+06}. 

Positive systems have gained increasing attention in the control
literature during the last decade. See for example 
\cite{Willems76,Farina+00,Kaczorek00}. Feedback stabilization of 
positive linear systems was studied in \cite{DeLeenheer+01}.
Stabilizing static output feedback controllers were
parameterized using linear programming in \cite{Rami+07,Rami11} and
extensions to input-ouput gain optimization were given in \cite{Briat13}.
Tanaka and Langbort
\cite{Tanaka+11} proved that the input-output gain of positive systems
can be evaluated using a diagonal quadratic storage
function and utilized this for $H_{\infty}$ optimization of
decentralized controllers in terms of semi-definite programming. A
related contribution is \cite{Najson13}, that proved a discrete
time Kalman-Yakubovich-Popov (KYP) lemma for positive
systems, with a different proof. 

The paper is structured as follows:
Section~\ref{sec:notation} introduces notation. Stability
criteria for positive systems are cited in
section~\ref{sec:stability}. These results are not new, but stated on
a form convenient for later use and explained with emphasis on scalability.
Section~\ref{sec:Hinf} extends the stability results to input-output
performance. The analysis results are then exploited in
section~\ref{sec:stabilization} for synthesis of
stabilizing and optimal controllers using distributed linear programming. Section~\ref{sec:LFD} extends the techniques to
positively dominated transfer functions. Section~\ref{sec:PQP}
explains how Lyapunov inequalities for positive systems can be
verified using methods that scale linearly with the number of states
and interconnections. Similar methods are used in section~\ref{sec:KYP} to prove
a more general version of the KYP lemma for positive systems.
The paper ends with conclusions and bibliography.

\section{Notation}
\label{sec:notation}
Let $\realR_+$ denote the set of nonnegative real numbers. For
$x\in\realR^n$, let $|x|\in\realR^n_+$ be the element-wise abolute
value. The notation ${\bf 1}$ denotes a column vector with all entries
equal to one.
The inequality \hbox{$X>0$ ($X\ge0$)} 
means that all elements of the matrix (or vector) $X$ are positive
(nonnegative). For a symmetric matrix $X$, the inequality $X\succ0$
means that the matrix is positive definite. The matrix
$A\in\realR^{n\times n}$ is said to be \emph{Hurwitz} if all
eigenvalues have negative real part. It is \emph{Schur} if all
eigenvalues are strictly inside the unit circle. Finally, the matrix
is said to be \emph{Metzler} if all off-diagonal elements are
nonnegative. 
The notation $\RHinfty$ represents the set of rational functions with
real coefficients and without poles in the closed right half plane. The set
of $n\times m$ matrices with elements in $\RHinfty$ is denoted $\RHinfty^{n\times m}$.

\section{Distributed Stability Verification}
\label{sec:stability}
The following well known characterizations of stability will be used extensively:
\begin{prp}
  Given a Metzler matrix $A\in\realR^{n\times n}$, the following statements are equivalent:
  \begin{description}
  \item[(\theprp.1)] The matrix $A$ is Hurwitz.
  \item[(\theprp.2)] There exists a ${\xi}\in\realR^{n}$ such that
    ${\xi}>0$ and $A{\xi}<0$.
  \item[(\theprp.3)] There exists a $z\in\realR^{n}$ such that $z>0$
    and $z^TA<0$.
  \item[(\theprp.4)] There exists a \emph{diagonal} matrix $P\succ0$ such that    
    $A^TP+PA\prec 0$.
  \item[(\theprp.5)] The matrix $-A^{-1}$ exists and has nonnegative entries.
  \end{description}
Moreover, if ${\xi}=({\xi}_1,\ldots,{\xi}_n)$ and $z=(z_1,\ldots,z_n)$ satisfy the
conditions of (\theprp.2) and (\theprp.3) respectively, then
$P=\diag(z_1/{\xi}_1,\ldots,z_n/{\xi}_n)$ satisfies the conditions of (\theprp.4).
\label{prp:contstab}
\end{prp}

\begin{figure}
  \begin{center}
    \includegraphics[width=0.25\hsize]{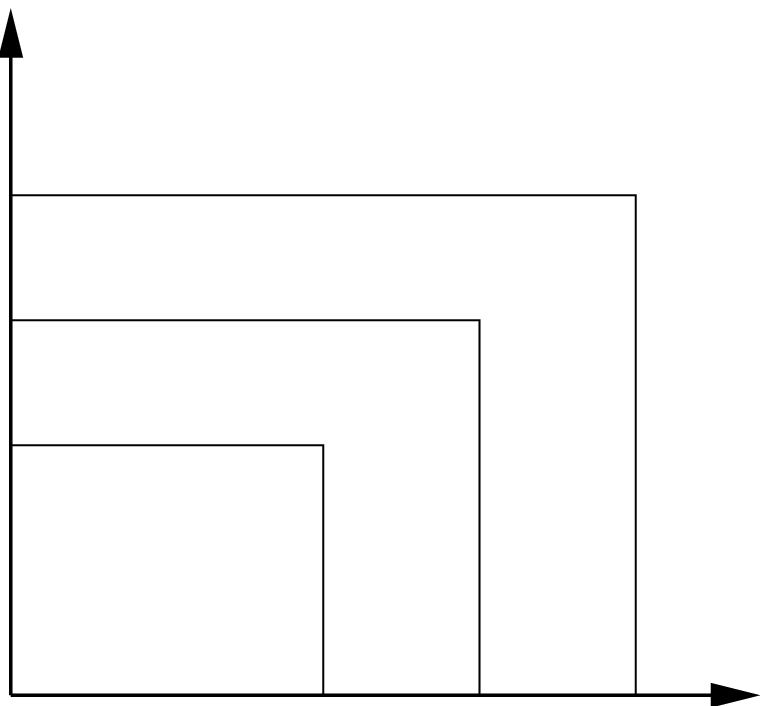}
    \hfill\includegraphics[width=0.25\hsize]{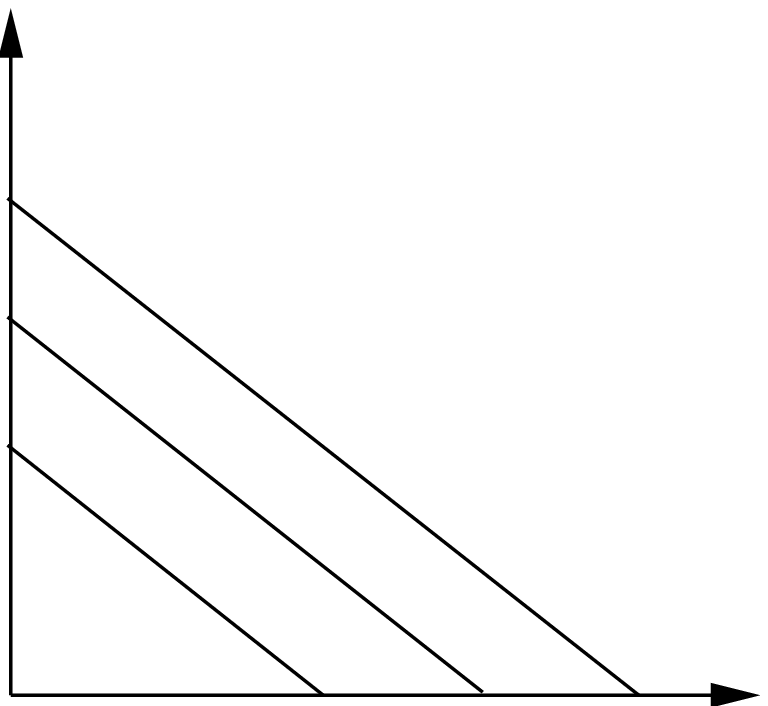}
    \hfill\includegraphics[width=0.25\hsize]{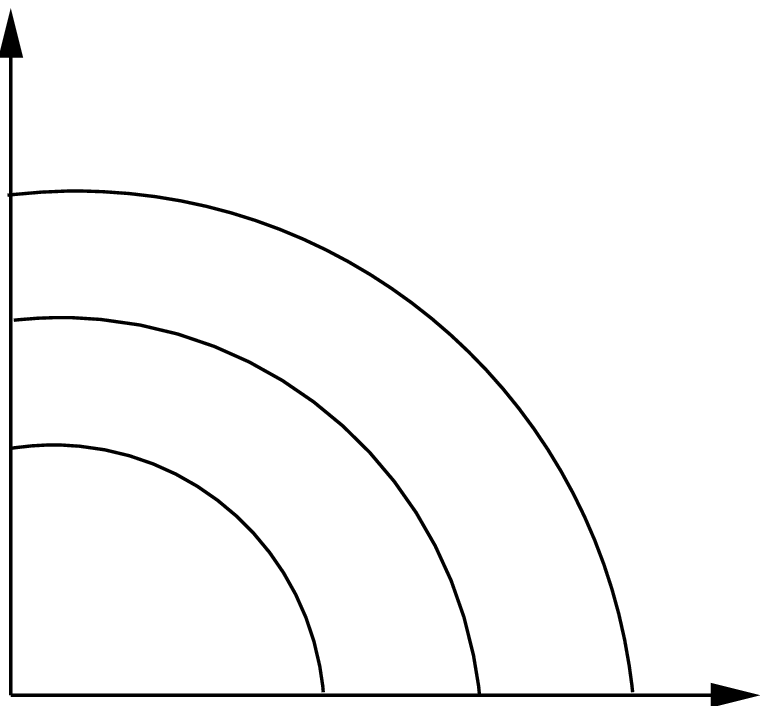}
  \end{center}
  \caption{Level curves of Lyapunov functions corresponding to the
    conditions (\ref{prp:contstab}.2), (\ref{prp:contstab}.3) and (\ref{prp:contstab}.4) in
    Proposition~\ref{prp:contstab}:\break If $A\xi<0$, then
    $V(x)=\max_i(x_i/\xi_i)$ is a
    Lyapunov function with rectangular level curves. If $z^TA<0$, then
    $V(x)=z^Tx$ is a linear 
    Lyapunov function. Finally if $A^TP+PA\prec 0$ and $P\succ0$, then
    $V(x)=x^TPx$ is a quadratic Lyapunov function for the system $\dot{x}=Ax$.}
  \label{fig:level}
\end{figure}
\begin{rmk}
  Each of the conditions (\ref{prp:contstab}.2), (\ref{prp:contstab}.3) and (\ref{prp:contstab}.4) corresponds to a
  Lyapunov function of a specific form. See Figure~\ref{fig:level}.
\end{rmk}

\begin{rmk}
  One of the main observations of this paper is that verification and
  synthesis of positive control systems can be done with methods that
  scale linearly with the number of interconnections. For stability,
  this claim follows directly from Proposition~\ref{prp:contstab}:
  Given $\xi$, verification of the inequality $A\xi<0$ requires a
  number of scalar additions and multiplications that is directly
  proportional to the number of nonzero elements in the matrix $A$. In
  fact, the search for a feasible $\xi$ also scales linearly, since
  integration of the differential equation $\dot{\xi}=A\xi$ with
  $\xi(0)=\xi_0$ for an arbitrary $\xi_0>0$ generates a feasible
  $\xi(t)$ in finite time provided that $A$ is Metzler and Hurwitz.
\end{rmk}

\begin{pf*}{Proof of Proposition~\ref{prp:contstab}.}
  The equivalence between (\theprp.1), (\theprp.2), (\theprp.4) and (\theprp.5) is the
equivalence between the statements $G_{20}$, $I_{27}$, $H_{24}$ and $N_{38}$ in
\cite[Theorem~6.2.3]{Berman+94}. The equivalence between (\ref{prp:contstab}.1) and
(\ref{prp:contstab}.3) is obtained by applying the equivalence between (\ref{prp:contstab}.1) and
(\ref{prp:contstab}.2) to the transpose of $A$. Moreover, if
${\xi}=({\xi}_1,\ldots,{\xi}_n)$ and $z=(z_1,\ldots,z_n)$ satisfy the
conditions of (\ref{prp:contstab}.2) and (\ref{prp:contstab}.3) respectively, then 
$P=\diag(z_1/{\xi}_1,\ldots,z_n/{\xi}_n)$ gives $(A^TP+PA)\xi=A^Tz+PA\xi<0$ so the symmetric matrix $A^TP+PA$ is Hurwitz and (\ref{prp:contstab}.4) follows.
\end{pf*}

\begin{ex*}{Linear transportation network.}
\label{ex:buffers}
  \begin{figure}
    \begin{center}
      \psfrag{V0}[l]{$x_1$}
      \psfrag{V1}[l]{$x_2$}
      \psfrag{V2}[l]{$x_3$}
      \psfrag{V3}[l]{$x_4$}
      \psfrag{w}[l]{$w$}
      \includegraphics[width=.7\hsize]{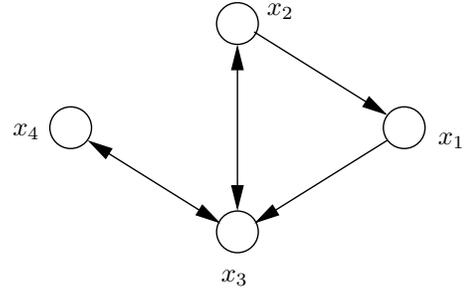}
    \end{center}
    \caption{A graph of an interconnected system. In Example~1 the
      interpretation is a transportation network and each arrow
      indicates a transportation link. In Example~2 the
      interpretation is instead a vehicle formation and each arrow
      indicates the use of a distance measurement.}
    \label{fig:vehicles}
  \end{figure}
Consider a dynamical system interconnected according to the graph illustrated in Figure~\ref{fig:vehicles}:
{\footnotesize\begin{align}
  \begin{bmatrix}
    \dot{x}_1\\\dot{x}_2\\\dot{x}_3\\\dot{x}_4
  \end{bmatrix}
  =
  \begin{bmatrix}
    -1-\ell_{31}&\ell_{12}&0&0\\
    0&-\ell_{12}-\ell_{32}&\ell_{23}&0\\
    \ell_{31}&\ell_{32}&-\ell_{23}-\ell_{43}&\ell_{34}\\
    0&0&\ell_{43}&-4-\ell_{34}
  \end{bmatrix}
  \begin{bmatrix}
    x_1\\x_2\\x_3\\x_4
  \end{bmatrix}
\label{eqn:buffers}
\end{align}
}The model could for example be used to describe a transportation
network connecting four buffers. The states $x_1,x_2,x_3,x_4$ represent the
contents of the buffers and the parameter $\ell_{ij}$ determines the
rate of transfer from buffer $j$ to buffer $i$. Such transfer
is necessary to stabilize the content of the second and third buffer.

Notice that the dynamics has the form $\dot{x}=Ax$
where $A$ is a Metzler matrix provided that every $\ell_{ij}$ is
nonnegative. Hence, by  
Proposition~\ref{prp:contstab}, stability is equivalent to existence
of numbers ${\xi}_1,\ldots,{\xi}_4>0$ such that
{\footnotesize\begin{align*}
  \begin{bmatrix}
    -1-\ell_{31}&\ell_{12}&0&0\\
    0&-\ell_{12}-\ell_{32}&\ell_{23}&0\\
    \ell_{31}&\ell_{32}&-\ell_{23}-\ell_{43}&\ell_{34}\\
    0&0&\ell_{43}&-4-\ell_{34}
  \end{bmatrix}
  \begin{bmatrix}
    {\xi}_1\\{\xi}_2\\{\xi}_3\\{\xi}_4
  \end{bmatrix}
< \begin{bmatrix}
    0\\0\\0\\0
  \end{bmatrix}
\end{align*}
}Given these numbers, stability can be verified by a distributed test where 
the first buffer verifies the first inequality,
the second buffer verifies the second and so on. In particular, the
relevant test for each buffer only involves parameter values at the
local node and the neighboring nodes, so a global model is not needed anywhere.
\end{ex*}

\begin{ex*}{Vehicle formation (or distributed Kalman filter).}
\label{ex:vehicles}
Another system structure, which can be viewed as a dual of the
previous one, is the following:
\begin{align}
  \begin{cases}
    \dot{x}_1=-x_1+\ell_{13}(x_3-x_1)\\
    \dot{x}_2=\ell_{21}(x_1-x_2)+\ell_{23}(x_3-x_2)\\  
    \dot{x}_3=\ell_{32}(x_2-x_3)+\ell_{34}(x_4-x_3)\\  
    \dot{x}_4=-4x_4+\ell_{43}(x_3-x_4)  
  \end{cases}
\label{eqn:vehicles}
\end{align}
This model could for example be used to describe a formation of four
vehicles. The parameters $\ell_{ij}$ represent position
adjustments based on distance measurements between the vehicles.
The terms $-x_1$ and $-4x_4$ reflect that the first and fourth vehicle
can maintain stable positions on their own, but the second and third vehicle rely
on the distance measurements for stabilization. Again,
stability can be verified by a distributed test where the first
vehicle verifies the first inequality, the second vehicle verifies the
second inequality and so on. 
\end{ex*}
A discrete time counterpart to Proposition~\ref{prp:contstab} is given
next:
\begin{prp}
  For $B\in\realR_+^{n\times n}$, the following statements are equivalent:
\begin{description}
  \item[(\theprp.1)] The matrix $B$ is Schur stable.
  \item[(\theprp.2)] There is a ${\xi}\in\realR^{n}$ such that
    ${\xi}>0$ and $B{\xi}<{\xi}$.
  \item[(\theprp.3)] There exists a $z\in\realR^{n}$ such that
    $z>0$ and $B^Tz<z$.
  \item[(\theprp.4)] There is a \emph{diagonal} $P\succ0$ such that
    $B^TPB\prec P$. 
  \item[(\theprp.5)] $(I-B)^{-1}$ exists and has nonnegative entries.
\end{description}
Moreover, if ${\xi}=({\xi}_1,\ldots,{\xi}_n)$ and $z=(z_1,\ldots,z_n)$
satisfy the
conditions of (\theprp.2) and (\theprp.3) respectively, then
$P=\diag(z_1/{\xi}_1,\ldots,z_n/{\xi}_n)$ satisfies the conditions of (\theprp.4).
\label{prp:discretestab}
\end{prp}

\begin{pf}
The equivalence between (\ref{prp:discretestab}.1) and
(\ref{prp:discretestab}.5) is proved by \cite[Lemma~6.2.1]{Berman+94}.
Setting $A=B-I$ gives the equivalence between
(\ref{prp:discretestab}.2), (\ref{prp:discretestab}.3) and
(\ref{prp:discretestab}.5) from the equivalence between
(\ref{prp:contstab}.2), (\ref{prp:contstab}.3) and
(\ref{prp:contstab}.5). 

Suppose ${\xi}=({\xi}_1,\ldots,{\xi}_n)$ and
$z=(z_1,\ldots,z_n)$ satisfy the conditions of
(\ref{prp:discretestab}.2) and (\ref{prp:discretestab}.3)
respectively. Set $P=\diag(z_1/{\xi}_1,\ldots,z_n/{\xi}_n)$. Then
  \begin{align*}
    B^TPB{\xi}&<B^TP{\xi}=B^Tz<P^{-1/2}z<P{\xi}
  \end{align*}
  so $B^TPB-P$ is Hurwitz and (\ref{prp:discretestab}.4) follows.
Finally, (\ref{prp:discretestab}.4) shows that $x^TPx$ is a positive
definite Lyapunov function for the system $x^+=Bx$, so
(\ref{prp:discretestab}.1) follows from (\ref{prp:discretestab}.4).
\end{pf}

\section{Input-Output Performance of Positive Systems}
\label{sec:Hinf}

We will now move beyond stability and discuss input-output
performance using induced norms. Given $M\in\realR^{r\times m}$, define the induced matrix norm
\begin{align*}
\|M\|_{p-\ind}&=\sup_{w\in\realR^m\setminus\{0\}}\frac{|Mw|_p}{|w|_p}
\end{align*}
where $|w|_p=(|w_1|^p+\cdots+|w_m|^p)^{1/p}$.
Assuming that $M$ has nonnegative entries we have
\begin{align*}
  \|M\|_{1-\ind}&<\gamma&&\text{ if and only if }&  M^T{\bf 1}&<\gamma {\bf 1}\\
  \|M\|_{\infty-\ind}&<\gamma&&\text{ if and only if }&  M{\bf 1}&<\gamma {\bf 1}
\end{align*}
For an $r\times m$ transfer matrix $\mathbf{G}(s)=C(sI-A)^{-1}B+D$, let $g(t)=Ce^{At}B+D\delta(t)$ be the corresponding impulse response. With $w\in\mathbf{L}_p^m[0,\infty)$, let $g*w\in\mathbf{L}_p^r[0,\infty)$ be
the convolution of $g$ and $w$ and define the induced norms
\begin{align*}
  \|g\|_{p-\ind}&=\sup_{w\in\mathbf{L}_p^m[0,\infty)}\frac{\|g*w\|_p}{\|w\|_p}
\end{align*}
where $\|w\|_p=\left(\sum_k\int_0^\infty|w_k(t)|^pdt\right)^{1/p}$.
A remarkable feature of positive
systems is that induced norms are determined by the static gain:

\begin{thm}
Let $g(t)=Ce^{At}B+D\delta(t)$ where $Ce^{At}B\ge0$ for $t\ge0$ and $D\ge0$, while
  $A$ is Hurwitz. Then $\|g\|_{p-\ind}=\|\mathbf{G}(0)\|_{p-\ind}$ for
  $p=1$, $p=2$ and $p=\infty$. In particular, if $g$ is scalar, then
  $\|g\|_{p-\ind}=\mathbf{G}(0)$ for all $p\in[1,\infty]$.
\label{thm:induced}
\end{thm}

\begin{pf}
It is well known that
$\|g\|_{2-\ind}=\max_{\omega}\|\mathbf{G}(i\omega)\|_{2-\ind}$ for
general linear time-invariant systems. When $g(t)\ge0$, the maximum must be attained at $\omega=0$ since 
\begin{align*}
  |\mathbf{G}({i\omega})w|&\le\int_0^\infty\Big|g(t)e^{-i\omega t}\Big|dt\cdot|w|\\
  &=\int_0^\infty g(t)dt\cdot|w|=\mathbf{G}(0)|w|
\end{align*}
for every $w\in\complexC^m$. This completes the proof for $p=2$. For $p=1$, the fact follows from
the calculations
  \begin{align*}
    \|y\|_1
    &=\sum_k\int_0^\infty\bigg|\sum_l\int_0^tg_{kl}(t-\tau)w_l(\tau)d\tau\bigg|_1dt\\
    &\le\sum_k\int_0^\infty\sum_l\int_0^tg_{kl}(t-\tau)|w_l(\tau)|d\tau dt\\
    &=\sum_{k,l}\int_0^\infty\bigg(\int_{\tau}^\infty
    g_{kl}(t-\tau)dt\bigg)|w_l(\tau)\big|d\tau\\
    &=\sum_{k,l}\left(\int_0^\infty g_{kl}(t)dt\right)\|w_l\|_1\\
    &=\sum_{k,l}\mathbf{G}_{kl}(0)\|w_l\|_1\\
    &\le\max_l\left(\sum_k\mathbf{G}_{kl}(0)\right)\|w\|_1\\
    &=\|\mathbf{G}(0)\|_{1-\ind}\cdot\|w\|_1
  \end{align*}
  with equality when
  $\|\mathbf{G}(0)\|_{1-\ind}\cdot\|w\|_1=\|\mathbf{G}(0)w\|_1$.
  Similarly, for $p=\infty$, 
  \begin{align*}
    \|y\|_\infty&=\max_{k,t}\left|\sum_l\int_0^\infty g_{kl}(\tau)w_l(t-\tau)d\tau\right|\\
    &\le\max_{k}\left(\sum_l\int_0^\infty g_{kl}(\tau)d\tau\right)\|w\|_\infty\\
    &=\max_{k}\left(\sum_l\mathbf{G}_{kl}(0)\right)\|w\|_\infty\\
    &=\|\mathbf{G}(0)\|_{\infty-\ind}\cdot\|w\|_\infty
  \end{align*}
  with equality when $w_l(t)$ has the same value for all $l$ and $t$. Hence the desired equality 
  \begin{align*}
    \|g\|_{p-\ind}&=\|\mathbf{G}(0)\|_{p-\ind}
  \end{align*}
  has been proved for $p=1$, $p=2$ and $p=\infty$. In particular, if $g$ is scalar, then
  \begin{align}
    \|g\|_{p-\ind}&=\mathbf{G}(0).
  \label{eqn:p-ind}
  \end{align}
  The Riesz-Thorin
  convexity theorem \cite[Theorem~7.1.12]{Hormander85} shows that
  $\|g\|_{p-\ind}$ is a convex function of $p$ for $1\le p\le\infty$,
  so (\ref{eqn:p-ind}) must hold for all $p\in[1,\infty]$.
\end{pf}

State-space conditions for input-output performance will now be established in parallel to the previous stability 
conditions:
\begin{thm}
  Let $g(t)=Ce^{At}B+D\delta(t)$ where
  $A\in\realR^{n\times n}$ is Metzler and $B\in\realR_+^{n\times
    m}$, $C\in\realR_+^{r\times n}$, $D\in\realR_+^{r\times m}$. Then the
  following statements are equivalent: 
  \begin{description}
  \item[(\theprp.1)] The matrix $A$ is Hurwitz and $\|g\|_{\infty-\ind}<\gamma$.\\[-2mm]
  \item[(\theprp.2)] There exists $\xi\in\realR_+^n$ such that 
  \begin{align}
  \begin{bmatrix}
    A&B\\C&D
  \end{bmatrix}\begin{bmatrix}\xi\\{\bf 1}\end{bmatrix}
  <\begin{bmatrix}0\\\gamma {\bf 1}\end{bmatrix}.
  \label{eqn:xi}
  \end{align}
  \end{description}
Moreover, if $\xi$ satisfies (\ref{eqn:xi}), then $-\xi<x(t)<\xi$ for all solutions to the
equation $\dot{x}=Ax+Bw$ with $x(0)=0$ and $\|w\|_\infty\le1$.
\label{thm:L1}
\end{thm}
\begin{pf*}{Proof of Theorem~\ref{thm:L1}.}
  $A$ is Metzler, so $e^{At}\ge0$ and the assumptions of
  Theorem~\ref{thm:induced} hold. Hence
  $\|g\|_{\infty-\ind}<\gamma$ can equivalently be written $\|D-CA^{-1}B\|_{\infty-\ind}<\gamma$ or
  \begin{align}
  (D-CA^{-1}B){\bf 1} < \gamma {\bf 1}.
  \label{eqn:ABCD}
  \end{align}
  Assume that (\ref{thm:L1}.2) holds. Then $A$ is Hurwitz by Proposition~\ref{prp:contstab}. Multiplying the inequality $A\xi+B{\bf 1}<0$ with the non-positive matrix $CA^{-1}$ from the left gives $C\xi+CA^{-1}B{\bf 1}\ge0$. Subtracting this from the inequality $C\xi+D{\bf 1}<\gamma{\bf 1}$ gives (\ref{eqn:ABCD}), so (\ref{thm:L1}.1) follows. 
  
  Conversely, suppose that (\ref{thm:L1}.1) and therefore (\ref{eqn:ABCD}) holds. By Proposition~\ref{prp:contstab} there exists $x>0$ such that $Ax<0$. Define $\xi=x-A^{-1}B$.
    Then $\xi\ge x>0$. Moreover
    \begin{align*}
      A\xi+B=Ax&<0
    \end{align*}
    If $x$ is sufficiently small, we also get $C\xi+D{\bf 1}<\gamma{\bf 1}$
    so (\ref{thm:L1}.2) follows.

   To prove the last statement, suppose that $\xi$ satisfies
   (\ref{eqn:xi}) and define $x$, $y$ and $z$ by
   \begin{align}
     \dot{y}&=Ay+u& y(0)&=-\xi\label{eqn:y}\\
     \dot{x}&=Ax+Bw& x(0)&=0 \label{eqn:x}\\
     \dot{z}&=Az+v& z(0)&=\xi \label{eqn:z}
   \end{align}
   where $\|w\|_\infty\le1$, $u=A\xi$ and $v=-A\xi$. Then the
   solutions of (\ref{eqn:y}) and (\ref{eqn:z}) are constantly equal
   to $-\xi$ and $\xi$ respectively. Moreover, the inequalities
   \begin{align*}
     u\le Bw \le v
   \end{align*}
   follow from (\ref{eqn:xi}). Together with the assumption that $A$ is Metzler, gives that
   $y(t)\le x(t)\le z(t)$ for all $t$. This completes the proof.
\end{pf*}

\begin{thm}
  Suppose that $g(t)=Ce^{At}B+D\delta(t)$ where
  $A\in\realR^{n\times n}$ is Metzler and $B\in\realR_+^{n\times
    m}$, $C\in\realR_+^{r\times n}$, $D\in\realR_+^{r\times m}$. Then the
  following statements are equivalent: 
  \begin{description}
  \item[(\theprp.1)] The matrix $A$ is Hurwitz and $\|g\|_{1-\ind}<\gamma$.\\[-2mm]
  \item[(\theprp.2)] There exists $p\in\realR_+^n$ such that 
  \begin{align}
  \begin{bmatrix}
    A&B\\C&D
  \end{bmatrix}^T\begin{bmatrix}p\\{\bf 1}\end{bmatrix}
  <\begin{bmatrix}0\\\gamma {\bf 1}\end{bmatrix}.\label{eqn:p}
  \end{align}
  \end{description}
  Moreover, if $p$ satisfies (\ref{eqn:p}), then all solutions to the
  equation $\dot{x}=Ax+Bw$ with $x(0)=0$ satisfy
  \begin{align}
    p^T|x(t)|+\int_0^t|Cx+Dw|d\tau&\le \gamma \int_0^t|w|d\tau
  \label{eqn:L1dissip}
  \end{align}
  with equality only if $w$ is identically zero.
  \label{thm:L1induced}
\end{thm}
\begin{rmk}
  The first part of Theorem~\ref{thm:L1} and Theorem~\ref{thm:L1induced}
  previously appeared in \cite{Briat13}.
\end{rmk}

\begin{pf}
  By Theorem~\ref{thm:induced}, the inequality $\|g\|_{1-\ind}<\gamma$ can
  equivalently be written $\|D-CA^{-1}B\|_{1-\ind}<\gamma$ or 
  \begin{align}
  (D-CA^{-1}B)^T{\bf 1} < \gamma {\bf 1}.
  \label{eqn:DCBA}
  \end{align}
Assume that (\ref{thm:L1induced}.1) holds. By Proposition~\ref{prp:contstab} there exists $z>0$ such that $z^TA<0$. Define $p=z-A^{-T}C^T$.
    Then $p\ge z>0$. Moreover
    \begin{align*}
      A^Tp+C^T=A^Tz&<0
    \end{align*}
    If $z$ is sufficiently small, we also get $B^Tp+D^T{\bf 1}<\gamma{\bf 1}$
    so (\ref{thm:L1induced}.2) follows.

    Conversely, suppose that (\ref{thm:L1induced}.2) holds. Then $A$ is Hurwitz by Proposition~\ref{prp:contstab}. Consider
    any solutions to
    \begin{align*}
       \dot{x}&=Ax+Bw&x(0)=0\\
       \dot{y}&=Ay+B|w|&y(0)=0.
    \end{align*}
    $A$ is Metzler, so $|x(t)|\le y(t)$ for all $t\ge0$. Multiplying
    the transpose of (\ref{eqn:p}) by $(y,|w|)$ from the right gives
    \begin{align*}
       p^T\dot{y}+Cy+D|w|\ge \gamma|w|.
    \end{align*}
    Integrating of $t$ and using that $|x(t)|\le y(t)$ 
    gives (\ref{eqn:L1dissip}). Then (\ref{thm:L1induced}.1) follows
    as $t\to\infty$ and the proof is complete.
\end{pf}

A discrete time counterpart of Theorem~\ref{thm:L1} and
Theorem~\ref{thm:L1induced} is given without proof:

\begin{thm}
  Given matrices $A,B,C,D\ge0$, let
  \begin{align*}
    g(t)&=\left\{
      \begin{array}{ll}
        \!\!D&t=0\\
        \!\!CA^{t-1}B&t=1,2,\ldots 
      \end{array}\right.
  \end{align*}
Then the
  following two statements are equivalent: 
  \begin{description}
  \item[(\theprp.1)] The matrix $A$ is Schur and $\|g\|_{\infty-\ind}<\gamma$.\\[-2mm]
  \item[(\theprp.2)] There exists $\xi\in\realR_+^n$ such that 
  \begin{align}
  \begin{bmatrix}
    A&B\\C&D
  \end{bmatrix}\begin{bmatrix}\xi\\{\bf 1}\end{bmatrix}
  <\begin{bmatrix}\xi\\\gamma{\bf 1}\end{bmatrix}.
  \label{eqn:dxi}
  \end{align}
  \end{description}
If $\xi$ satisfies (\ref{eqn:dxi}), then $-\xi<x(t)<\xi$ for all solutions to the
equation $x(t+1)=Ax(t)+Bw(t)$ with $x(0)=0$ and $\|w\|_\infty\le1$.

  The following two statements are also equivalent: 
  \begin{description}
  \item[(\theprp.3)] The matrix $A$ is Schur and $\|g\|_{1-\ind}<\gamma$.\\[-2mm]
  \item[(\theprp.4)] There exists $p\in\realR_+^n$ such that 
  \begin{align}
  \begin{bmatrix}
    A&B\\C&D
  \end{bmatrix}^T\begin{bmatrix}p\\{\bf 1}\end{bmatrix}
  <\begin{bmatrix}p\\\gamma {\bf 1}\end{bmatrix}.\label{eqn:dp}
  \end{align}
  \end{description}
  Moreover, if $p$ satisfies (\ref{eqn:dp}), then all solutions to the
  equation $x(t+1)=Ax(t)+Bw(t)$ with $x(0)=0$ satisfy
  \begin{align*}
    p^T|x(t)|+\sum_{\tau=0}^t|Cx(\tau)+Dw(\tau)|&\le \gamma \sum_{\tau=0}^t|w(\tau)|
  \end{align*}
  with equality only if $w$ is identically zero.
  \label{thm:discreteperf}
\end{thm}

\section{Distributed Control Synthesis by Linear Programming}
\label{sec:stabilization}

Equipped with scalable analysis methods for stability and performance,
we are now ready to consider synthesis of controllers by
distributed optimization. We will start
by re-visiting an example of section~\ref{sec:stability}. 

\begin{ex}
Consider again the transportation network (\ref{eqn:buffers}), this
time with the flow parameters $\ell_{31}=2$, $\ell_{34}=1$ and $\ell_{43}=2$ fixed:
\begin{align}
  \begin{bmatrix}
    \dot{x}_1\\\dot{x}_2\\\dot{x}_3\\\dot{x}_4
  \end{bmatrix}
  =
  \begin{bmatrix}
    -3&\ell_{12}&0&0\\
    0&-\ell_{12}-\ell_{32}&\ell_{23}&0\\
    2&\ell_{32}&-\ell_{23}-2&1\\
    0&0&2&-5
  \end{bmatrix}
  \begin{bmatrix}
    x_1\\x_2\\x_3\\x_4
  \end{bmatrix}
\label{eqn:transportcontrol}
\end{align}
We will ask the question how to find the remaining parameters $\ell_{12}$, $\ell_{23}$ and $\ell_{32}$ in the interval
$[0,1]$ such that the closed loop system (\ref{eqn:transportcontrol})
becomes stable. According to Proposition~\ref{prp:contstab}, stability
is equivalent to existence of $\xi_1,\ldots,\xi_4>0$ such that
\begin{align*}
  \begin{bmatrix}
    -3&\ell_{12}&0&0\\
    0&-\ell_{12}-\ell_{32}&\ell_{23}&0\\
    2&\ell_{32}&-\ell_{23}-2&1\\
    0&0&2&-5
  \end{bmatrix}
  \begin{bmatrix}
    \xi_1\\\xi_2\\\xi_3\\\xi_4
  \end{bmatrix}&<0
\end{align*}
At first sight, this looks like a difficult problem due to
multiplications between the two categories of parameters. However, a
closer look suggests the introduction of new variables:
$\mu_{12}:=\ell_{12}\xi_2$, $\mu_{32}:=\ell_{32}\xi_2$ and
$\mu_{23}:=\ell_{23}\xi_3$. The problem then reduces to linear
programming: Find $\xi_1,\xi_2,\xi_3,\xi_4>0$ and
$\mu_{12},\mu_{32},\mu_{23}\ge0$ such that
\begin{align*}
  &\begin{bmatrix}
    -3&0&0&0\\
    0&0&0&0\\
    2&0&-2&1\\
    0&0&2&-5
  \end{bmatrix}
  \begin{bmatrix}
    \xi_1\\\xi_2\\\xi_3\\\xi_4
  \end{bmatrix}
  +
  \begin{bmatrix}
     1&0&0\\
    -1&-1&1\\
     0&1&-1\\
     0&0&0
  \end{bmatrix}\begin{bmatrix}\mu_{12}\\\mu_{32}\\\mu_{23}\end{bmatrix}
<0\\
&\mu_{12}\le\xi_2\qquad \mu_{32}\le\xi_2\qquad \mu_{23}\le\xi_3
\end{align*}
with the solution $(\xi_1,\xi_2,\xi_3,\xi_4)=(0.5,0.5,1.69,0.87)$ and
$(\mu_{12},\mu_{32},\mu_{23})=(0.5,0.5,0)$. The corresponding stabilizing gains can
then be computed as
\begin{align*}
  \ell_{12}&=\mu_{12}/\xi_2=1&
  \ell_{32}&=\mu_{32}/\xi_2=1&
  \ell_{23}&=\mu_{23}/\xi_3=0
\end{align*}
\label{ex:transynt}
\end{ex}
The idea can be generalized into the following theorem:

\begin{thm}
Let $\mathcal{D}$ be the set of $m\times m$ diagonal matrices with
entries in $[0,1]$. Suppose that ${A}+{E}L{F}$ is Metzler and
${C}+{G}L{F}\ge0$, $B+ELH \ge0$, $D+GLH \ge0$ for all $L\in\mathcal{D}$.
Let $g_L(t)$ be the impulse response of
\begin{align*}
  (C+GLF)[sI-(A+ELF)]^{-1}(B+ELH)+D+GLH
\end{align*}
If ${F}\ge0$, then the following two conditions are equivalent: 
  \begin{description}
  \item[(\theprp.1)] There exists $L\in\mathcal{D}$ with
  ${A}+{E}L{F}$ is Hurwitz and $\|g_L\|_{\infty-\ind}<\gamma$.
  \item[(\theprp.2)] There exist $\xi\in\realR_+^n$, $\mu\in\realR_+^m$ with
    \begin{align*}
      {A}{\xi}+{B}{\bf 1}+{E}{\mu}&<0\\
      {C}{\xi}+{D}{\bf 1}+{G}{\mu}&<\gamma{\bf 1}\\
      {F}{\xi}+{H}{\bf 1}&\ge {\mu}
    \end{align*}
  \end{description}
  Moreover, if $\xi,\mu$ satisfy (\theprp.2), then (\theprp.1) holds for every $L$
  such that $\mu=LF{\xi}+L{H}{\bf 1}$.
\label{thm:Hinfsyn}
\end{thm}

\begin{figure}
\begin{center}
\psfrag{v1}[c]{$x$}
\psfrag{v2}[c]{$w$}
\psfrag{z1}[c]{$\dot{x}$}
\psfrag{z2}[c]{$z$}
\psfrag{y}[c]{$y$}
\psfrag{u}[c]{$u$}
\psfrag{P}[c]{$\begin{bmatrix}A&B&E\\C&D&G\\F&H&0\end{bmatrix}$}
\psfrag{C}[c]{$L$}
\psfrag{s}[c]{$s^{-1}$}
\includegraphics[width=.9\hsize]{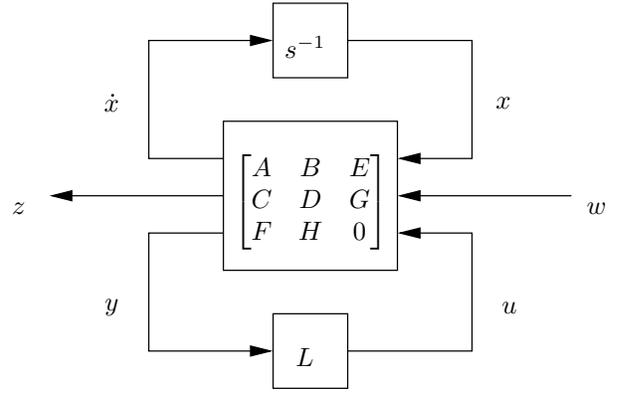}
\end{center}
\caption{Theorem~\ref{thm:Hinfsyn} shows how to determine the
  diagonal matrix $L$ that minimizes the $L_{\infty}$-induced gain from $w$ to
$z$. Theorem~\ref{thm:Hinfsyndual}
shows how to minimize the
$L_{1}$-induced gain. If $w$ and $z$ are scalar, both gains are equal
to the the standard $H_{\infty}$-norm. Extension to the case of system
matrix with non-zero lower right corner is straightforward, but omitted from this
paper.} 
\label{fig:LFT}
\end{figure}
\begin{rmk}
  If the diagonal elements of
  $\mathcal{D}$ are restricted to $\realR_+$ instead of $[0,1]$, then the
  condition $ {F}{\xi}+{H}{\bf 1}\ge {\mu}$ is replaced by $ {F}{\xi}+{H}{\bf 1}\ge0$. 
\end{rmk}

\begin{rmk}
  When the matrices have a sparsity pattern corresponding to a graph,
  each row of the vector inequalities in (\theprp.2) can be verified
  separately to get a distributed performance test. 

  Also finding a solution to the linear programming problem
  can be done with distributed methods, where each node in the graph
  runs a local algorithm involving only local variables and
  information exchange only with its neighbors. For example, given a
  stable Metzler matrix $A$, consider
  the problem to find a stability certificate $\xi>0$ satisfying
  $A\xi<0$. This can be done in a distributed way by simulating the system using Euler's
  method until the state is close to a dominating
  eigenvector of the $A$. Then it must satisfy the conditions on $\xi$.
\end{rmk}

\begin{rmk}
It is interesting to compare our results with the analysis and
synthesis methods proposed by Tanaka and Langbort in \cite{Tanaka+11}
and Briat in \cite{Briat13}. Our
mathematical treatment has much in common with theirs. However,
none of them is discussing scalable design, nor verification, of
distributed controllers. Moreover, our ``static output feedback''
expression $A+ELF$ is significantly more general than the ``state feedback''
expression $A+BL$ used in both those references. This gives us a
higher degree of flexibility, particularly in the specification of
distributed controllers. On the other hand, their parametrization has
the advantage that the Metzler property of the closed loop system
matrix can be enforced as a constraint in the synthesis procedure,
rather than being verified a priori for all $L\in\mathcal{D}$. 
\end{rmk}

\begin{pf}
    Suppose (\ref{thm:Hinfsyn}.1) holds. Then, according to
    Theorem~\ref{thm:L1}, there exists $\xi\in\realR_+^n$ such that 
  \begin{align}
  \begin{bmatrix}
    A+ELF&B+ELH\\C+GLF&D+GLH
  \end{bmatrix}\begin{bmatrix}\xi\\{\bf 1}\end{bmatrix}
  <\begin{bmatrix}0\\\gamma {\bf 1}\end{bmatrix}.
  \label{eqn:AELF}
  \end{align}
  Setting ${\mu}=LF{\xi}+L{H}{\bf 1}$ gives (\ref{thm:Hinfsyn}.2). Conversely, suppose that (\ref{thm:Hinfsyn}.2) holds. Choose $L\in\mathcal{D}$ to get
  ${\mu}=LF{\xi}+L{H}{\bf 1}$. Then (\ref{eqn:AELF}) holds and
  (\ref{thm:Hinfsyn}.1) follows by Theorem~\ref{thm:L1}.
\end{pf}
Theorem~\ref{thm:Hinfsyn} was inspired by the transportation network
in Example~3, where non-negativity of $F$ is
natural assumption. However, this condition would fail in a vehicle formation
problem, where control is based on distance measurements. For such
problems, the following dual formulation is useful:
\begin{thm}
Let $\mathcal{D}$ be the set of $m\times m$ diagonal matrices with
entries in $[0,1]$. Suppose that ${A}+{E}L{F}$ is Metzler and
${C}+{G}L{F}\ge0$, $B+ELH \ge0$, $D+GLH \ge0$ for all $L\in\mathcal{D}$.
Let $g_L(t)$ be the impulse response of
\begin{align*}
  (C+GLF)[sI-(A+ELF)]^{-1}(B+ELH)+D+GLH
\end{align*}
If the matrices ${B},{D}$ and ${E}$ have
nonnegative coefficients, then the following two conditions are equivalent: 
  \begin{description}
  \item[(\theprp.1)] There exists $L\in\mathcal{D}$ with
  ${A}+{E}L{F}$ is Hurwitz and $\|g_L\|_{1-\ind}<\gamma$.
  \item[(\theprp.2)] There exist $p\in\realR_+^n$, $q\in\realR_+^m$ with
    \begin{align*}
      {A}^T{p}+{C}^T{\bf 1}+{F}^T{q}&<0\\
      {B}^T{p}+{D}^T{\bf 1}+{H}^T{q}&<\gamma{\bf 1}\\
      {E}^T{p}+{G}^T {\bf 1}&\ge{q}
    \end{align*}
  \end{description}
  Moreover, if $p,q$ satisfy (\theprp.2), then (\theprp.1) holds for every $L$
  such that $q=LE^T{p}+L{G}^T {\bf 1}$.
\label{thm:Hinfsyndual}
\end{thm}
\begin{pf}
  The proof is analogous to the proof of Theorem~\ref{thm:Hinfsyn}.
\end{pf}
\begin{ex*}{Disturbance rejection in vehicle formation.}
  Consider the vehicle formation model
\begin{align}
  \begin{cases}
    \dot{x}_1=-x_1+\ell_{13}(x_3-x_1)+w\\
    \dot{x}_2=\ell_{21}(x_1-x_2)+\ell_{23}(x_3-x_2)+w\\  
    \dot{x}_3=\ell_{32}(x_2-x_3)+\ell_{34}(x_4-x_3)+w\\  
    \dot{x}_4=-4x_4+\ell_{43}(x_3-x_4)+w
  \end{cases}
\label{eqn:vehiclepert}
\end{align}
where $w$ is an external disturbance acting on the vehicles. 
Our problem is to find feedback gains
  gains $\ell_{ij}\in[0,1]$ that stabilize the
  formation and minimize the gain from $w$ to $x$. The problem can
  be solved by applying Theorem~\ref{thm:Hinfsyndual} with 
  {\small\begin{align*}
    A&=\diag\{-1,0,0,-4\}&
    C&=
    \left(\begin{array}{rrrr}
  \!\! 1& 1& 1& 1
    \end{array}\right)\\[2mm] 
    E&=
    \left(\begin{array}{rrrrrr}
  \!\!1& 0& 0& 0& 0& 0\\
       0&1&1& 0& 0& 0\\
       0& 0& 0&1&1& 0\\
       0& 0& 0& 0& 0&1
    \end{array}\right)&
    D&=0\\[2mm]
    L&=\diag\{\ell_{13},\ell_{21},\ell_{23},\ell_{32},\ell_{34},\ell_{43}\}\\[2mm]
    F&=
    \left(\begin{array}{rrrr}
  \!\! -1& 0& 1& 0\\
        1&-1& 0& 0\\
        0&-1& 1& 0\\
        0& 1&-1& 0\\
        0& 0&-1& 1\\
        0& 0& 1&-1
    \end{array}\right)&
    H&=
    \left(\begin{array}{r}
  \!\! 0\\
        0\\
        0\\
        0\\
        0\\
        0
    \end{array}\right)
  \end{align*}
}Solutions for three different cases are illustrated in Figure~\ref{fig:Enrico}.
\begin{figure}
  \begin{center}
      \includegraphics[width=.3\hsize]{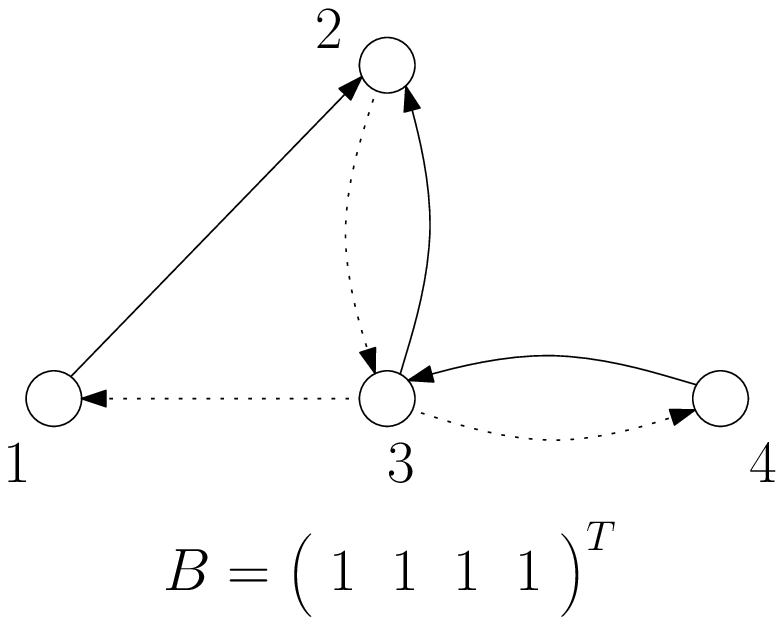}\quad 
      \includegraphics[width=.3\hsize]{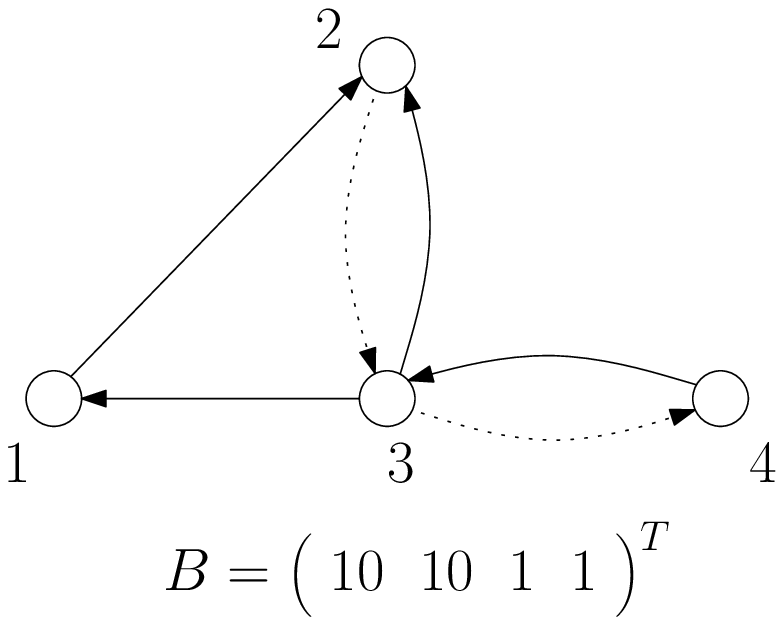}\quad 
      \includegraphics[width=.3\hsize]{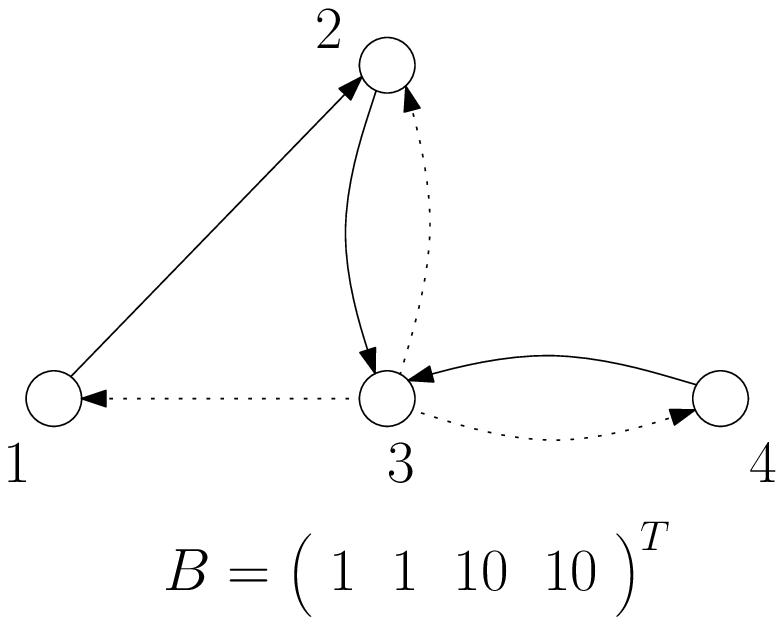}
  \end{center}
  \caption{Illustration of optimal gains for disturbance rejection in
    a vehicle formation. When $B=(1\;1\;1\;1)^T$, all four vehicles face unit disturbances and
    the optimal $L=\diag\{0,1,1,0,1,0\}$ illustrated by arrows in the
    left diagram gives $\gamma= 4.125$. Apparently, the first vehicle
    should ignore
  the distance to the third vehicle, while the third vehicle should
  ignore the second vehicle and the fourth should ignore the third. The middle
diagram illustrates a situation where the disturbances on vehicle 1
and 2 are ten times bigger. Then the minimal value $\gamma= 15.562$ is attained
with $L=\diag\{1,1,1,0,1,0\}$, so the first vehicle should use
distance measurements to the third. The converse situation in the right
diagram gives $\gamma=12.750$ for $L=\diag\{0,1,0,1,1,0\}$.}
  \label{fig:Enrico}
\end{figure}
\end{ex*}

\section{Positively Dominated Systems}
\label{sec:LFD}

So far, the emphasis has been on state space models. However, for many
applications input-output models are more natural as a starting point.
In this section, we will therefore extend the main ideas of the previous
sections to such models. First we need to define a notion of
positivity for input-output models. One option would be to work with
non-negative impulse responses like in Theorem~\ref{thm:induced}.
However, to verify for a given rational transfer function that the
impulse response is non-negative has proved to be NP-hard! See \cite{Blondel02} for the discrete time problem and
\cite{Bell+10} for continuous time. Instead we will use the
following definition.

$\mathbf{G}\in\RHinfty^{m\times n}$ is called
\emph{positively dominated} if every matrix entry satisfies 
$|\mathbf{G}_{jk}(i\omega)|\le \mathbf{G}_{jk}(0)$ for all
$\omega\in\realR$. 
The set of all such matrices is denoted 
$\DHinfty^{m\times n}$. The essential scalar frequency inequality can be
tested by semi-definite programming, since 
$|b(i\omega)/a(i\omega)|\le b(0)/a(0)$ holds for $\omega\in\realR$
if and only if the polynomial
$|a(i\omega)|^2b(0)^2-|b(i\omega)|^2a(0)^2$ can be written as a sum of squares.

Some properties of positively dominated transfer functions follow immediately:
\begin{prp}
  Let $\mathbf{G},\mathbf{H}\in\DHinfty^{n\times n}$. Then $\mathbf{G}\mathbf{H}\in\DHinfty^{n\times n}$ and 
  $a\mathbf{G}+b\mathbf{H}\in\DHinfty^{n\times n}$ when $a,b\in\realR_+$. 
  Moreover $\|\mathbf{G}\|_{\infty}=\|\mathbf{G}(0)\|$.
\label{prp:lfd}
\end{prp}

The following property is also fundamental:
\begin{thm}
  Let $\mathbf{G}\in\DHinfty^{n\times n}$.
  Then $(I-\mathbf{G})^{-1}\in\DHinfty^{n\times n}$ if and
  only if $\mathbf{G}(0)$ is Schur. 
\label{thm:lfdstab}
\end{thm}

\begin{pf}
  That $(I-\mathbf{G})^{-1}$ is stable and positively dominated
  implies that $[I-\mathbf{G}(0)]^{-1}$ exists and is nonnegative, so
  $\mathbf{G}(0)$ must be Schur according to
  Proposition~\ref{prp:discretestab}. On the other hand, if
  $\mathbf{G}(0)$ is Schur 
  we may choose ${\xi}\in\realR_+$ and $\epsilon>0$ with $\mathbf{G}(0){\xi}<(1-\epsilon){\xi}$. Then
  for every $z\in\complexC^n$ with $0<|z|<\xi$ and $s\in\complexC$
  with $\hbox{Re }s\ge0$ we have
  \begin{align*}
    |\mathbf{G}(s)^tz|\le \mathbf{G}(0)^t|z|<(1-\epsilon)^t|z|&&\hbox{for }t=1,2,3,\ldots
  \end{align*}
  Hence $\sum_{k=0}^\infty \mathbf{G}(s)^tz$ is convergent and bounded
  above by $\sum_{k=0}^\infty \mathbf{G}(0)^t|z|=[I-\mathbf{G}(0)]^{-1}|z|$. The sum of the
  series solves the equation $[I-\mathbf{G}(s)]\sum_{k=0}^\infty \mathbf{G}(s)^tz=z$, so
  therefore $\sum_{k=0}^\infty
  \mathbf{G}(s)^tz=[I-\mathbf{G}(s)]^{-1}z$. This proves
  $(I-\mathbf{G})^{-1}$ is stable and positively dominated and the
  proof is complete.
\end{pf}

\begin{figure}
\begin{center}
\psfrag{v1}[c]{$x$}
\psfrag{v2}[c]{$w$}
\psfrag{z1}[c]{}
\psfrag{z2}[c]{$z$}
\psfrag{y}[c]{$y$}
\psfrag{u}[c]{$u$}
\psfrag{P}[c]{$\begin{bmatrix}\mathbf{A}&\mathbf{B}&\mathbf{E}\\\mathbf{C}&\mathbf{D}&0\\\mathbf{F}&0&0\end{bmatrix}$}
\psfrag{C}[c]{$L$}
\includegraphics[width=.8\hsize]{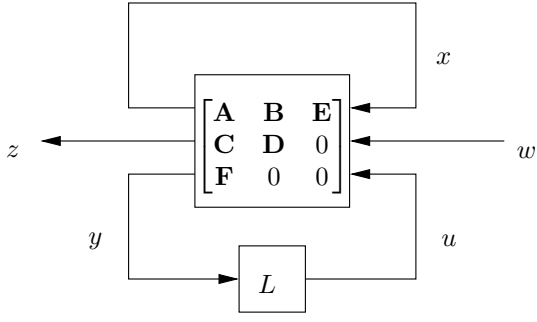}
\end{center}
\caption{Theorem~\ref{thm:lfdgain}
shows how to determine the
  diagonal matrix $L$ that minimizes the $L_1$-induced gain from $w$ to
$z$. Extension to the case of a matrix with all non-zero blocks
is straightforward, but omitted from this paper. The same is true for
the transpose version corresponding to Theorem~\ref{thm:Hinfsyn}.} 
\label{fig:LFT2}
\end{figure}

Theorem~\ref{thm:Hinfsyndual} has the following counterpart for
positively dominated systems, as illustrated in Figure~\ref{fig:LFT2}.
\begin{thm}
Let $\mathcal{D}$ be the set of $m\times m$ diagonal matrices with
entries in $[0,1]$, while $\mathbf{B}\in\DHinfty^{n\times k}$, $\mathbf{C}\in\DHinfty^{l\times n}$,
  $\mathbf{D}\in\DHinfty^{l\times k}$,
  $\mathbf{E}\in\DHinfty^{n\times m}$ and
  $\mathbf{F}\in\RHinfty^{m\times n}$. Suppose
$\mathbf{A}+\mathbf{E}L\mathbf{F}\in\DHinfty^{n\times
  n}$ for all $L\in\mathcal{D}$. \\[-3mm]

Then the following two conditions are equivalent: 
  \begin{description}
  \item[(\theprp.1)] 
  There is $L\in\mathcal{D}$ with $(I-\mathbf{A}-\mathbf{E}{L}\mathbf{F})^{-1}\in\DHinfty^{n\times n}$
  and
  $\|\mathbf{C}(I-\mathbf{A}-\mathbf{E}L\mathbf{F})^{-1}\mathbf{B}+\mathbf{D}\|_{1-\ind}<\gamma$.
  \item[(\theprp.2)] 
%
There exist $p\in\realR_+^n$, $q\in\realR_+^m$ with
    \begin{align*}
      \mathbf{A}(0)^Tp+\mathbf{C}(0)^T {\bf 1}+\mathbf{F}(0)^Tq&<p\\
      \mathbf{B}(0)^Tp+\mathbf{D}(0)^T {\bf 1}&<\gamma{\bf 1}\\
      \mathbf{E}(0)^Tp&\ge q
    \end{align*}
\end{description}
  If $p,q$ satisfy (\theprp.2), then (\theprp.1) holds for every $L$
  such that $q=L\mathbf{E}(0)^Tp$.
\label{thm:lfdgain}
\end{thm}

\begin{pf}
  Proposition~\ref{prp:lfd} and Theorem~\ref{thm:lfdstab} show that (\ref{thm:lfdgain}.1) holds if and only if
  $\mathbf{A}(0)-\mathbf{E}(0){L}\mathbf{F}(0)$ is Schur and
  \begin{align*}
    \|\mathbf{C}[I-\mathbf{A}(0)-\mathbf{E}(0)L\mathbf{F}(0)]^{-1}\mathbf{B}(0)+\mathbf{D}(0)\|_{1-\ind}<\gamma
  \end{align*}
  According to Theorem~\ref{thm:discreteperf}, this is true if and
  only if there exists $p\in\realR_+^n$ such that 
  \begin{align*}
    \begin{bmatrix}
      \mathbf{A}(0)+\mathbf{E}(0)L\mathbf{F}(0)&\mathbf{B}(0)\\
      \mathbf{C}(0)&\mathbf{D}(0)
    \end{bmatrix}^T\begin{bmatrix}p\\{\bf 1}\end{bmatrix}
    <\begin{bmatrix}p\\\gamma{\bf 1}\end{bmatrix}.
  \end{align*}
  This is equivalent to (\ref{thm:lfdgain}.2) if we set 
  $q=L\mathbf{E}(0)^T{p}$, so the desired equivalence between
  (\ref{thm:lfdgain}.1) and (\ref{thm:lfdgain}.2) follows. 
\end{pf}

\begin{ex*}{Formation of vehicles with inertia.}
  In Example~3, the inputs and disturbances were supposed to have an
  immediate impact on the vehicle velocities, i.e. the inertia of the
  vehicles was neglected. Alternatively, a model that takes the
  inertia into account can be stated as follows:
  \begin{align*}
    \ddot{x}_i&=\sum_{j}{}\ell_{ij}(x_j-x_i)+u_i+w_i&
    i=1,\ldots,N
  \end{align*}
  where $u_i$ is a control force, $w_i$ is a
  disturbance force and ${\ell}_{ij}$ is the spring constant between the
  vehicles $i$ and $j$. Suppose that local control laws
  $u_i=-k_ix_i-d_i\dot{x}$ are given and consider
  the problem to 
  find spring constants ${\ell}_{ij}\in[0,\overline{\ell}_{ij}]$ that minimize the gain from
  $w_1$ to $x_1$.

  The closed loop system has the frequency domain description
  \begin{align*}
    &\bigg(s^2+d_is+k_i+\sum_j\overline{\ell}_{ij}\bigg)X_i(s)\\
    &=\sum_j\bigg({\ell}_{ij}X_j(s)+(\overline{\ell}_{ij}-{\ell}_{ij})X_i(s)\bigg)+W_i(s).
  \end{align*}
  Similarly to Example~3, we write this on matrix form as
  \begin{align*}
    X&=(\mathbf{A}+\mathbf{E}L\mathbf{F})X+\mathbf{B}W
  \end{align*}
  The transfer matrices $\mathbf{B}$, $\mathbf{E}$ and
  $\mathbf{A}+\mathbf{E}L\mathbf{F}$ are
  positively dominated for all $L\in\mathcal{D}$ provided that $d_i\ge k_i+\sum_j\overline{\ell}_{ij}$. 
  Hence Theorem~\ref{thm:lfdgain} can then be applied to find the
  optimal spring constants. Notice that $\ell_{ij}$ and $\ell_{ji}$
  must be optimized separately, even though by symmetry they must be equal at optimum.
\end{ex*}

\section{Scalable Verification of the Lyapunov Inequality}
\label{sec:PQP}

In the preceding sections we have derived scalable conditions for
verification of stability and optimality, using generalizations of the linear
inequalities in (\ref{prp:contstab}.2) and (\ref{prp:contstab}.3) of
Proposition~\ref{prp:contstab}. To address multi-variable systems using linear programming, the
natural performance measures have been input-output gains with
signals measured $L_1$-norm or $L_\infty$-norm.

A more well-known alternative, used in the classical $H_{\infty}$ control theory, is to measure signals with $L_2$-norm. This
was done in \cite{Tanaka+11} using generalizations of condition
(\ref{prp:contstab}.4), however without discussion of scalability aspects.
The purpose of the next theorem is to show that for positive systems also verification of semi-definite
inequalities, like the Lyapunov inequality
$A^TP+PA\prec0$, can be decomposed into tests that scale linearly with the number of
non-zero matrix entries. 
\begin{thm}
  A symmetric Metzler matrix with $m$ non-zero entries above the
  diagonal is negative semi-definite if and only if
  it can be written as a sum of $m$ negative semi-definite matrices,
  each of which has only four non-zero entries.
\label{thm:posdecomp}
\end{thm}

The proof of Theorem~\ref{thm:posdecomp} will be based on the
following minor modification of \cite[Theorem~3.1]{Kim+03}:

\begin{prp}[Positive Quadratic Programming]
  Suppose
  $M_0,\ldots,M_K$ are Metzler and
  $b_1,\ldots,b_K\in\realR$. Then
  \begin{align}
    \begin{array}{llcll}
      \max&x^TM_0x&=&\max&\trace(M_0X)\\[1mm]
      x\in\realR^n_+&x^TM_kx\,\ge\,b_k&&
      X\succeq0&\trace(M_kX)\ge b_k\\
      &k=1,\ldots,K&&&k=1,\ldots,K
    \end{array}
  \label{eqn:PQP}
  \end{align}
  The value on the right hand side remains the same if the condition
  $X\succeq0$ is relaxed to $X\in\mathbb{X}$, where $\mathbb{X}$ is the set of symmetric
  matrices $(x_{ij})\in\realR^{n\times n}$ satisfying $x_{ii}\ge0$ and
  $x_{ij}^2\le x_{ii}x_{jj}$ for all $i,j$. 
  Moreover, if there exists a matrix $X$ in the interior of $\mathbb{X}$ with
  $\trace(M_kX)\ge b_k$ for every $k$, then the maximum of
  (\ref{eqn:PQP}) is 
  equal to the minimum of $-\sum_k\tau_kb_k$ over
  $\tau_1,\ldots,\tau_K\ge0$ such that $M_0+\sum_k\tau_kM_k\preceq0$. 
\label{prp:PQP}
\end{prp}

\begin{rmk}
  The problem on the right of (\ref{eqn:PQP}) is always convex and readily solvable by
  semidefinite programming. The problem on the left is generally not a
  convex program, since the matrices $M_k$ may be indefinite. However,
  the maximization on the left is concave in
  $(x_1^2,\ldots,x_n^2)$ \cite{MegretskiPC10}. This is because every product $x_ix_j$ is the
  geometric mean of two such variables, hence concave \cite[p. 74]{Boyd+04}.
\end{rmk}

\begin{rmk}
  The second statement of Proposition~\ref{prp:PQP} is important for scalability, since the
  condition $X\in\mathbb{X}$ has a natural decomposition and only 
  entries of $X$ that correspond to non-zero entries of $M_k$ need to
  be taken into account. 
\end{rmk}

\begin{pf*}{Proof of Proposition~\ref{prp:PQP}.}
Every $x$ satisfying the constraints on the left hand side of (\ref{eqn:PQP}) corresponds
to a matrix $X=xx^T$ satisfying the constraints on the right hand
side. This shows that the right hand side of (\ref{eqn:PQP}) is at least as big as the left. 

On the other hand, let $X=(x_{ij})$ be a positive definite matrix. In particular, the diagonal elements
$x_{11},\ldots,x_{nn}$ are non-negative and
$x_{ij}\le\sqrt{x_{ii}x_{jj}}$.  
Let $x=(\sqrt{x_{11}},\ldots,\sqrt{x_{nn}})$. Then the matrix $xx^T$
has the same diagonal elements as $X$, but has off-diagonal elements
$\sqrt{x_{ii}x_{jj}}$ instead of $x_{ij}$. The fact that $xx^T$ has
off-diagonal elements at least as big as those of $X$, together with
the assumption that the matrices $M_k$ are Metzler, gives
$x^TM_kx\ge\trace(M_kX)$ for $k=1,\ldots,K$. This shows that the left
hand side of (\ref{eqn:PQP}) is at least as big as the right. Nothing
changes if $X$ is not positive definite but $X\in\mathbb{X}$, so the
second statement is also proved.

For the last statement, note that the conditions $\trace(M_kX)\ge b_k$
are linear in $X$, so strong duality holds
\cite[Theorem~28.2]{Rockafellar97} and the
right hand side of (\ref{eqn:PQP}) 
has a finite maximum if and only if
  $M_0+\sum_{k=1}^K\tau_kM_k\preceq0$ for some $\tau_1,\ldots,\tau_K\ge0$.
\end{pf*}

\begin{pf*}{Proof of Theorem~\ref{thm:posdecomp}.}
  Let $\mathcal{E}$ be the set of indices $(k,l)$ of non-zero
  off-diagonal entries in $M$. Define
  \begin{align*}
    \mathbb{X}_{\mathcal{E}}&=\left\{X\in\realR^{n\times n}:\,\,\hbox{$[e_k\;\;e_l]^TX[e_k\;\;e_l]\succeq0$
      for all }(k,l)\in\mathcal{E}\right\}
  \end{align*}
  where $e_1,\ldots,e_n$ are the unit vectors in $\realR^n$.
  If $M$ is negative semi-definite, then
  {\small\begin{align*}
    0&=\max_{|x|\le1}x^TMx=\max_{X\in\mathbb{X}_{\mathcal{E}}}\trace(MX)\\
    &=\min_{N_{kl}\succeq0}\max_{X\in\realR^{n\times n}}\trace(MX)+\!\!\!\!\sum_{(k,l)\in\mathcal{E}}\!\!\trace\left(N_{kl}[e_k\;\;e_l]^TX[e_k\;\;e_l]\right)\\
    &=\min_{N_{kl}\succeq0}\max_{X\in\realR^{n\times n}}\trace\left[\left(M+\sum_{(k,l)\in\mathcal{E}}[e_k\;\;e_l]N_{kl}[e_k\;\;e_l]^T\right)X\right]
  \end{align*}
  }where $N_{kl}\in\realR^{2\times 2}$ for every $k$ and $l$.
  In particular, there exists a choice of the matrices $N_{kl}$ that
  makes $M+\sum_{(k,l)\in\mathcal{E}}[e_k\;\;e_l]N_{kl}[e_k\;\;e_l]^T=0$. This
  completes the proof.
\end{pf*}

\section{The KYP Lemma for Positive Systems}
\label{sec:KYP}

Input-output gain is certainly not the only way to quantify the
performance of a linear time-invariant system. A more general class of
specifications known as Integral Quadratic Constraints
\cite{meg+ran97} can be tested using the Kalman-Yakubovich-Popov
lemma. It is therefore of interest to see that the corresponding
result of
\cite{Tanaka+11} for positive systems can be generalized the following
way:
\begin{thm}
Let $A\in\realR^{n\times n}$ be Metzler and Hurwitz, while
$B\in\realR_+^{n\times m}$ and the pair $(-A,B)$ is stabilizable. Suppose that all entries of
$Q\in\realR^{(n+m)\times(n+m)}$ are nonnegative, except for the last
$m$ diagonal elements. Then the following statements are equivalent:\\[-2mm]
\begin{description}
\item[(\theprp.1)]\quad For $\omega\in[0,\infty]$ is is true that
  \begin{align*}
    \begin{bmatrix}(i\omega I-A)^{-1}B\\I\end{bmatrix}^*Q
    \begin{bmatrix}(i\omega I-A)^{-1}B\\I\end{bmatrix}\preceq 0
  \end{align*}\\[-3mm]
\item[(\theprp.2) ]$\quad\begin{bmatrix}-A^{-1}B\\I\end{bmatrix}^TQ
    \begin{bmatrix}-A^{-1}B\\I\end{bmatrix}\preceq 0$.\\[3mm]
\item[(\theprp.3) ]\quad There exists a diagonal $P\succeq0$ such that
  \begin{align*}
    Q+\begin{bmatrix}
        A^TP+PA&PB\\
        B^TP&0
      \end{bmatrix}\preceq 0
  \end{align*}
\item[(\theprp.4) ]\quad There exist $x,p\in\realR^{n}_+$, $u\in\realR^{m}_+$
  with $Ax+Bu\le0$,
  $$Q \begin{bmatrix}x\\u\end{bmatrix}+\begin{bmatrix}A^T\\B^T\end{bmatrix}p\le0$$
\end{description}
Moreover, if all inequalities are replaced by strict ones, then the
equivalences hold even without the stabilizability assumption.
\label{thm:posKYP}
\end{thm}

\begin{rmk}
  For $A=-1$, $B=0$, $Q=$ {\tiny$\begin{bmatrix}0&1\\1&0\end{bmatrix}$},
  condition (\ref{thm:posKYP}.1) holds, but not (\ref{thm:posKYP}.3).
  This demonstrates that the stabilizability of $(-A,B)$ is essential. 
\end{rmk}
\begin{rmk}
  Our statement of the KYP lemma for continuous and discrete time
  positive systems extends earlier versions of
  \cite{Tanaka+11,Najson13} in several respects: Non-strict
  inequality, more general $Q$ and a fourth equivalent condition in
  terms of linear programming rather than semi-definite programming.
\end{rmk}
\begin{pf}
One at a time, we will prove the implications
(\ref{thm:posKYP}.1) $\Rightarrow$ (\ref{thm:posKYP}.2) $\Rightarrow$
(\ref{thm:posKYP}.3) $\Rightarrow$ (\ref{thm:posKYP}.1)
and (\ref{thm:posKYP}.2) $\Leftrightarrow$ (\ref{thm:posKYP}.4).
Putting $\omega=0$ immediately gives (\ref{thm:posKYP}.2) from
(\ref{thm:posKYP}.1). 

Assume that (\ref{thm:posKYP}.2) holds. The
matrix $-A^{-1}$ is nonnegative, so
{\scriptsize$\begin{bmatrix}x\\w\end{bmatrix}^T$}$Q\,${\scriptsize$ 
  \begin{bmatrix}x\\w\end{bmatrix}$}$\le0$ for all $x\in\realR^n_+$,
  $w\in\realR_+^m$ with
  \begin{align}
    x&\le -A^{-1}Bw
  \label{eqn:xw}
  \end{align}
  The inequality (\ref{eqn:xw}) follows (by multiplication with
  $-A^{-1}$ from the left) from the constraint $0\le Ax+Bw$, which can also be written
  $0\le A_ix+B_iw$ for $i=1,\ldots,n$, 
  where $A_i$ and $B_i$ denote the $i$:th rows of $A$ and $B$
  respectively. For non-negative $x$ and $w$, this is
  equivalent to 
  \begin{align}
    0&\le x_i(A_ix+B_iw)&
    i&=1,\ldots,n
  \label{eqn:xwi2}
  \end{align}
  Hence (\ref{thm:posKYP}.2) implies {\scriptsize$\begin{bmatrix}x\\w\end{bmatrix}^T$}$Q\,${\scriptsize$
  \begin{bmatrix}x\\w\end{bmatrix}$}$\le0$ for
  $x\in\realR^n_+$,
  $w\in\realR_+^m$ satisfying (\ref{eqn:xwi2}). 
  Proposition~\ref{prp:PQP} will next be used to verify existence of
  $\tau_1,\ldots,\tau_n\ge0$ such that the quadratic form 
  \begin{align*}
    \sigma(x,w)=\Biggl[\!\!\begin{array}{c}x\\w\end{array}\!\!\Biggr]^T
    Q\Biggl[\!\!\begin{array}{c}x\\w\end{array}\!\!\Biggr]+\sum_i\tau_ix_i(A_ix+B_iw)
  \end{align*}
  is negative semi-definite. However, the application of
  Proposition~\ref{prp:PQP} requires existence of a positive
  definite $X$ such that all diagonal elements of
  \begin{align*}
    \begin{bmatrix}
      A&B
    \end{bmatrix}X
    \begin{bmatrix}
      I\\0
    \end{bmatrix}
  \end{align*}
  are positive. The pair
  $(-A,B)$ is stabilizable, so there exists $K$ that make all
  eigenvalues of $A+BK$ unstable and therefore
  $(A+BK)Z+Z(A+BK)^T=I$ has a symmetric positive definite solution
  $Z$. Hence the desired $X$ can be constructed as
  \begin{align*}
    X=
    \begin{bmatrix}
      Z&ZK^T\\KZ&*
    \end{bmatrix}
  \end{align*}
  where the lower right corner is chosen big enough to make $X\succ0$.

  Define $P=\diag(\tau_1,\ldots,\tau_n)\succeq0$. Then $\sigma$ being
  negative definite means that
  \begin{align*}
    Q+\begin{bmatrix}
        A^TP+PA&PB\\
        B^TP&0
      \end{bmatrix}\preceq 0
  \end{align*}
  so (\ref{thm:posKYP}.3) follows. 

Assume that (\ref{thm:posKYP}.3) holds. Integrating $\sigma(x(t),w(t))$ over time gives
  \begin{align*}
    0&\ge\int_0^\infty
    \left(\Biggl[\!\!\begin{array}{c}x\\w\end{array}\!\!\Biggr]^T
    Q\Biggl[\!\!\begin{array}{c}x\\w\end{array}\!\!\Biggr]
     +x^TP(Ax+Bw)\right)dt
  \end{align*}
  For square integrable solutions to $\dot{x}=Ax+Bw$, $x(0)=0$ we get
  \begin{align*}
    0&\ge\int_0^\infty
    \left(\Biggl[\!\!\begin{array}{c}x\\w\end{array}\!\!\Biggr]^T
    Q\Biggl[\!\!\begin{array}{c}x\\w\end{array}\!\!\Biggr]
     +\frac{d}{dt}(x^TPx/2)\right)dt\\
    &=\int_0^\infty
    \Biggl[\!\!\begin{array}{c}x(t)\\w(t)\end{array}\!\!\Biggr]^T
    Q\Biggl[\!\!\begin{array}{c}x(t)\\w(t)\end{array}\!\!\Biggr]dt
  \end{align*}
  which in frequency domain implies (\ref{thm:posKYP}.1). Hence
(\ref{thm:posKYP}.1) $\Rightarrow$ (\ref{thm:posKYP}.2) $\Rightarrow$
(\ref{thm:posKYP}.3) $\Rightarrow$ (\ref{thm:posKYP}.1). 

Assuming again (\ref{thm:posKYP}.2) gives, by
Proposition~\ref{prp:contstab}, existence of $u\in\realR^{m}_+$ such that
\begin{align*}
  \left(\begin{bmatrix}-A^{-1}B\\I\end{bmatrix}^TQ
    \begin{bmatrix}-A^{-1}B\\I\end{bmatrix}\right)u\le0
\end{align*}
Setting $x=-A^{-1}u$ gives $x\in\realR^{n}_+$ and
\begin{align*}
  \begin{bmatrix}-A^{-1}&A^{-1}B\\0&-I\end{bmatrix}^TQ\begin{bmatrix}x\\u\end{bmatrix}\ge0
\end{align*}
due to the sign structure of $Q$. Let $\begin{bmatrix}p^T&q^T\end{bmatrix}^T$ be the column on the left
hand side. Multiplying with $\begin{bmatrix}A&B\\0&I\end{bmatrix}^T$
from the left gives
\begin{align*}
  \begin{bmatrix}A^T&0\\B^T&I\end{bmatrix}\begin{bmatrix}p\\q\end{bmatrix}
  =-Q\begin{bmatrix}x\\u\end{bmatrix}
\end{align*}
and (\ref{thm:posKYP}.4) follows. 

Finally, suppose that (\ref{thm:posKYP}.4) holds. Then $x\ge
-A^{-1}Bu$. Multiplying the main inequality from the left with
$\begin{bmatrix}-B^TA^{-T}&I\end{bmatrix}$ gives
\begin{align*}
  0&\ge\begin{bmatrix}-A^{-1}B\\I\end{bmatrix}^TQ\begin{bmatrix}x\\u\end{bmatrix}
  =\left(\begin{bmatrix}-A^{-1}B\\I\end{bmatrix}^TQ
    \begin{bmatrix}-A^{-1}B\\I\end{bmatrix}\right)u
\end{align*}
and (\ref{thm:posKYP}.2) follows.

For strict inequalities, the proofs that  (\ref{thm:posKYP}.2)
$\Leftrightarrow$ (\ref{thm:posKYP}.4) and 
(\ref{thm:posKYP}.3) $\Rightarrow$ (\ref{thm:posKYP}.1) $\Rightarrow$ (\ref{thm:posKYP}.2)
remain the same. Assuming that
(\ref{thm:posKYP}.2) holds with strict inequality, we get
\begin{align*}
  \begin{bmatrix}-A^{-1}B\\I\end{bmatrix}^*(Q+\epsilon I)
    \begin{bmatrix}-A^{-1}B\\I\end{bmatrix}\preceq 0
\end{align*}
for some scalar $\epsilon>0$. Hence, there exists a diagonal $P\succeq0$ such that
\begin{align*}
  Q+\epsilon I+\begin{bmatrix}
        A^TP+PA&PB\\
        B^TP&0
      \end{bmatrix}\preceq 0
\end{align*}
Adding a small multiple of the identity to $P$ gives $P\succ0$ such that
\begin{align*}
  Q+\begin{bmatrix}
        A^TP+PA&PB\\
        B^TP&0
      \end{bmatrix}\prec 0
\end{align*}
so also (\ref{thm:posKYP}.3) holds with strict inequality. Hence the
proof is complete. 
\end{pf}

An analogous discrete time result is stated here and proved in the appendix:

\begin{thm}
Let $A\in\realR_+^{n\times n}$ be Schur, while
$B\in\realR_+^{n\times m}$ and the pair $(A,B)$ is anti-stabilizable. Suppose that all entries of
$Q\in\realR^{(n+m)\times(n+m)}$ are nonnegative, except for the last
$m$ diagonal elements. Then the following statements are equivalent:\\
\begin{description}
\item[(\theprp.1)]\quad For $\omega\in[0,\infty]$ is is true that
  \begin{align*}
    \begin{bmatrix}(e^{i\omega}I-A)^{-1}B\\I\end{bmatrix}^*Q
    \begin{bmatrix}(e^{i\omega}I-A)^{-1}B\\I\end{bmatrix}\preceq 0
  \end{align*}\\[-3mm]
\item[(\theprp.2) ]$\quad\begin{bmatrix}(I-A)^{-1}B\\I\end{bmatrix}^*Q
    \begin{bmatrix}(I-A)^{-1}B\\I\end{bmatrix}\preceq 0$.\\[3mm]
\item[(\theprp.3) ]\quad There exists a diagonal $P\succeq0$ such that
  \begin{align*}
    Q+\begin{bmatrix}
        A^TPA-P&A^TPB\\
        B^TPA&B^TPB
      \end{bmatrix}\preceq 0
  \end{align*}
\item[(\theprp.4) ]\quad There are $x,p\in\realR^{n}_+$, $u\in\realR^{m}_+$
  with $x\ge Ax+Bu$,
  $$Q \begin{bmatrix}x\\u\end{bmatrix}+\begin{bmatrix}A^T-I\\B^T\end{bmatrix}p\le0$$
\end{description}
Moreover, if all inequalities are taken to be strict, then the
equivalences hold even without the  anti-stabilizability assumption.
\label{thm:discrKYP}
\end{thm}

\begin{pf}
  The theorem can be proved in analogy with the proof of Theorem~\ref{thm:posKYP}.
  Alternatively, it can be derived from Theorem~\ref{thm:posKYP} using
  a bilinear transformation in the following way:

  Instead of $e^{i\omega}$, one can parametrize the unit circle as
  $\frac{1+i\omega}{1-i\omega}$. Hence (\ref{thm:discrKYP}.1) is
  equivalent to saying that
  \begin{align*}
    \begin{bmatrix}x\\u\end{bmatrix}^*Q\begin{bmatrix}x\\u\end{bmatrix}\le0
  \end{align*}
  for all solutions $(\omega,x,u)$ to the equation $\left(\frac{1+i\omega}{1-i\omega}I-A\right)x=Bu$.
  Alternatively, introducing
  \begin{align*}
    \widehat{A}&=(A-I)(A+I)^{-1}\\
    \widehat{B}&=2(A+I)^{-1}B\\
    \widehat{x}&=x+Ax+Bu\\
    S&=\begin{bmatrix}(A+I)^{-1}&-(A+I)^{-1}B\\0&I\end{bmatrix}\\
    \widehat{Q}&=S^TQS
  \end{align*}
  the condition can be re-written as the statement that
  \begin{align*}
    \begin{bmatrix}\widehat{x}\\u\end{bmatrix}^*\widehat{Q}\begin{bmatrix}\widehat{x}\\u\end{bmatrix}\le0
  \end{align*}
  for all solutions $(\omega,\widehat{x},\widehat{u})$ to the
  equation $(i\omega I-\widehat{A})\widehat{x}=\widehat{B}u$.
  According to Theorem~\ref{thm:posKYP}, this is equivalent to
  validity of the inequality for $\omega=0$,
  i.e. (\ref{thm:discrKYP}.2). It is also equivalent to existence of a diagonal $P\succeq0$ such that
  \begin{align*}
    \widehat{Q}+\begin{bmatrix}
        \widehat{A}^TP+P\widehat{A}&P\widehat{B}\\
        \widehat{B}^TP&0
      \end{bmatrix}\preceq 0
  \end{align*}
  Multiplying by {\small$\begin{bmatrix}A+I&B\\0&I\end{bmatrix}$} from
  the right and its transpose from the left, the matrix inequality (after trivial
  manipulations) becomes
  \begin{align*}
    Q+2\begin{bmatrix}
        A^TPA-P&A^TPB\\
        B^TPA&B^TPB
      \end{bmatrix}\preceq 0
  \end{align*}
  Replacing $2P$ by $P$ gives equivalence to
  (\ref{thm:discrKYP}.3). Also by Theorem~\ref{thm:posKYP}, it is equivalent to existence of
  $\widehat{x},p\in\realR^{n}_+$, $\widehat{u}\in\realR^{m}_+$ such that
  \begin{align*}
    \widehat{Q}\begin{bmatrix}\widehat{x}\\u\end{bmatrix}
    +\begin{bmatrix}\widehat{A}^T\\\widehat{B}^T\end{bmatrix}p\le0
  \end{align*}
  Left multiplication by $S^{-T}$ and substitution
  $(x,u)=S(\widehat{x},u)$ gives equivalence to (\ref{thm:discrKYP}.4).
\end{pf}

As an application of the equivalence between (\ref{thm:posKYP}.1) and
(\ref{thm:posKYP}.2), we consider an example devoted to optimal power
flow in an electrical network, a time-varying version of a problem
considered in\cite{lavaei+11}:

\begin{ex*}{Optimal power flow in an electrical network.}
\begin{figure}
  \begin{center}
    \psfrag{V0}{$\!v_4$}
    \psfrag{I0}{$i_4$}
    \psfrag{V1}{$\!v_1$}
    \psfrag{I1}{$i_1$}
    \psfrag{V2}{$\!v_2$}
    \psfrag{I2}{$i_2$}
    \psfrag{V3}{$\!\!v_3$}
    \psfrag{I3}{$i_3$}
    \includegraphics[width=0.65\hsize]{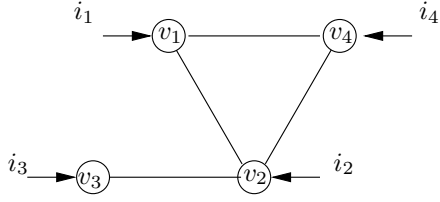}
  \end{center}
  \caption{Illustration of a dynamic power transmission network 
    with inductive transmission lines studied in Example~6.}
  \label{fig:powergraph}
\end{figure}
Consider a power transmission network as in
Figure~\ref{fig:powergraph}. The current from node $j$ to node $k$ is
governed by the voltage difference $v_j-v_k$ according to the
differential equation
\begin{align}
  L_{jk}\frac{di_{jk}}{dt}&=-R_{jk}i_{jk}+v_j(t)-v_k(t)  
\label{eqn:inductive}
\end{align}
and the external currents are determined by Kirchoff's law
\begin{align}
  \begin{cases}
  i_1(t)=-i_{41}(t)-i_{21}(t)\\
  i_2(t)=i_{21}(t)-i_{32}(t)-i_{42}(t)\\
  i_3(t)=i_{32}(t)\\
  i_4(t)=i_{41}(t)+i_{42}(t)
  \end{cases}
\label{eqn:Kirchoff}
\end{align}
The generation and consumption of power is subject to constraints of
the form
\begin{align}
  \frac{1}{T}\int_0^Ti_k(t)v_k(t)dt\le \overline{p}_k
\label{eqn:power}
\end{align}
If $k$ is a generator node, then $\overline{p}_k>0$ indicates
production capacity. Similarly, 
for loads $\overline{p}_k<0$ represents power demand. Transmission
lines have capacity constraints of the form
\begin{align}
  \frac{1}{T}\int_0^T|v_k(t)-v_j(t)|^2dt\le c_{kj}
\label{eqn:capacity}
\end{align}
Finally, the voltages are non-negative and subject to magnitude bounds
\begin{align}
  \underline{v}_k^2\le\frac{1}{T}\int_0^T v_k(t)^2dt\le\overline{v}_k^2
\label{eqn:magnitude}
\end{align}
We are now interested to minimize the resistive power losses in the
network subject to the given constraints:
\begin{align*}
  \begin{cases}
    \hbox{Minimize }\frac{1}{T}\sum_{k=1}^4\int_0^T i_k(t)v_k(t)dt\\[2mm]
    \hbox{subject to }(\ref{eqn:inductive})-(\ref{eqn:magnitude}) 
  \end{cases}
\end{align*}
Using the theory above, our goal is to prove that \emph{minimal losses can
be attained with constant voltages and currents.}

With line currents being states and voltage differences being inputs, this is a problem of the form
\begin{align}
  \begin{cases}
    \hbox{Maximize $\int_0^T$
    {\scriptsize$\begin{bmatrix}x\\u\end{bmatrix}^*$}
    $Q_0$
     {\scriptsize$\begin{bmatrix}x\\u\end{bmatrix}$}
    $dt$}\\[2mm]
    \hbox{subject to }\dot{x}(t)=Ax(t)+Bu(t)\\[2mm]
    \hbox{and $\int_0^T$
     {\scriptsize$\begin{bmatrix}x\\u\end{bmatrix}^*$}
    $Q_k$
     {\scriptsize$\begin{bmatrix}x\\u\end{bmatrix}$} 
    $dt\le q_k,\quad k=1,\ldots,m$ }
  \end{cases}
\label{eqn:dynopt}
\end{align}
where $A\in\realR^{n\times n}$ is Metzler and Hurwitz, 
$B\in\realR_+^{n\times m}$ and all entries of
$Q_k\in\realR^{(n+m)\times(n+m)}$ are nonnegative except possibly for
the last $m$ diagonal elements. To bring the problem on a form where
Theorem~\ref{thm:posKYP} can be applied, we will apply relaxation in
two different ways: The inequalities are handled using Lagrange
relaxation and the time interval is extended to $[0,\infty]$. This
brings the problem to the dual form
\begin{align*}
  \begin{cases}
  \hbox{Minimize }-\sum_{k=1}^m\tau_kq_k \hbox{ subject to
    $\tau_k\ge0$ such that}\\[2mm]
    \hbox{$\int_0^\infty$
    {\scriptsize$\begin{bmatrix}x\\u\end{bmatrix}^*$}
    $ (Q_0+\sum_k\tau_kQ_k)$
    {\scriptsize$\begin{bmatrix}x\\u\end{bmatrix}$}
    $dt\le0$}\\[2mm]
  \hbox{for all solutions to }\dot{x}(t)=Ax(t)+Bu(t).
  \end{cases}
\end{align*}
In frequency domain, this is written as
\begin{align*}
  \begin{cases}
  \hbox{Minimize }\sum_{k=1}^m\tau_kq_k \hbox{ subject to $\tau_k\ge0$ and}\\[2mm]
  \hbox{\small$\begin{bmatrix}(i\omega I-A)^{-1}B\\I\end{bmatrix}^*(Q_0+\sum_k\tau_kQ_k) 
    \begin{bmatrix}(i\omega I-A)^{-1}B\\I\end{bmatrix}\preceq 0$}
  \end{cases}
\end{align*}
The equivalence between (\ref{thm:posKYP}.1) and (\ref{thm:posKYP}.2)
in Theorem~\ref{thm:posKYP} shows that the bottleneck is always the frequency $\omega=0$, so
the problem takes the form
\begin{align*}
  &\min_{\tau_k\ge0}\sum_{k=1}^m\tau_kq_k \hbox{ subject to }
    M_0+\sum_k\tau_kM_k \preceq 0
\end{align*}
where
\begin{align*}
  M_k&=\begin{bmatrix}-A^{-1}B\\I\end{bmatrix}^*Q_k 
    \begin{bmatrix}-A^{-1}B\\I\end{bmatrix},&k&=1,\ldots,m
\end{align*}
are all Metzler. By Proposition~\ref{eqn:PQP} this can be restated as
\begin{align*}
  \max_u u^TM_0u \hbox{ subject to }u^TM_ku\le q_k
\end{align*}
or equivalently
\begin{align}
  \begin{cases}
  \hbox{Maximize 
  {\scriptsize$\begin{bmatrix}x\\u\end{bmatrix}^T$}
  $Q_0$
  {\scriptsize$\begin{bmatrix}x\\u\end{bmatrix}$}}\\[2mm]
  \hbox{subject to }Ax+Bu=0\\[2mm]
  \hbox{and 
  {\scriptsize$\begin{bmatrix}x\\u\end{bmatrix}^T$}
  $Q_k$
   {\scriptsize$\begin{bmatrix}x\\u\end{bmatrix}$}
  $\le q_k\quad k=1,\ldots,m$}
  \end{cases}
\label{eqn:statopt}
\end{align}
Notice that (\ref{eqn:statopt}) was obtained from (\ref{eqn:dynopt})
by relaxation, so the optimal value of (\ref{eqn:statopt}) must be at
least as good as the value of (\ref{eqn:dynopt}). At the same time,
(\ref{eqn:statopt}) is the special case of (\ref{eqn:dynopt})
obtained with constant values of the variables, so our goal has been achieved:
Minimal losses can be attained with constant voltages and currents.
\end{ex*}

\section{Conclusions}
The results above demonstrate that the monotonicity properties of
positive systems and positively dominated systems bring remarkable
benefits to control theory. Most important is
the opportunity for scalable verification and synthesis of
distributed control systems with optimal input-output performance. In particular, linear
programming solutions come with certificates that enable distributed
and scalable verification of global optimality, without access to a global model anywhere.

Many important problems remain open for future research. Here are two examples:
\begin{itemize}
\item How can the scalable methods for verification and synthesis be
  extended to monotone nonlinear systems?
\item How can the controller optimization be extended to scalable methods for optimization of dynamic
  controllers?
\end{itemize}

\section{Acknowledgments}
The author is grateful for suggestions and
comments by numerous colleagues, including the anonymous reviewers. In particular,
Theorem~\ref{thm:induced} was updated based on suggestions by Andrej
Ghulchak and Figure~\ref{fig:Enrico} was drawn by Enrico Lovisari. The work has been supported by the Swedish
Research Council through the Linnaeus Center LCCC.

\bibliography{regler,refs,rantzer}

\end{document}